\documentclass[a4paper, 12pt]{article}

\usepackage{amssymb,latexsym} \author{Johannes Sj\"ostrand,\\ \small
  CMLS, Ecole Polytechnique\\ \small FR-91128 Palaiseau, France \\
\small johannes@math.polytechnique.fr\\
\footnotesize UMR 7640 CNRS
}
\date{Dedicated to Hans Duistermaat}
\title{Eigenvalue distributions and Weyl laws for
  semi-classical non-self-adjoint operators in 2 dimensions}

\newtheorem{dref}{Definition}[section] 
\newtheorem{theo}[dref]{Theorem} 
\newtheorem{remark}[dref]{Remark}

\newenvironment{proof}{\par\noindent{{\bf Proof}}}{\hfill$\Box$
\medskip}

\newcommand{\ekv}[2]{\begin{equation}\label{#1}#2\end{equation}}

  \newcommand\iint{\int\hskip -2mm\int}

\newcommand{\no}[1]{(\ref{#1})} 
\begin{document}

\maketitle
\begin{abstract}   
In this note we compare two recent results about the distribution of
eigenvalues for semi-classical pseudodifferential operators in two
dimensions. For classes of analytic operators A.~Melin and the author
\cite{MeSj2} obtained a complex Bohr-Sommerfeld rule, showing that the
eigenvalues are situated on a distorted lattice. On the other hand,
with M.~Hager \cite{HaSj} we showed in any dimension that Weyl
asymptotics holds with probability close to 1 for small random
perturbations of the operator. In both cases the eigenvalues
distribute to leading order according two smooth densities and we show
here that the two densities are in general different.
\medskip
\par \centerline{\bf R\'esum\'e}
Dans cette note nous comparons deux r\'esultats r\'ecents sur la
distribution des valeurs propres d'op\'erateurs pseudodiff\'erentiels
en deux dimensions. Pour des classes d'op\'erateurs analytiques
A.~Melin et l'auteur \cite{MeSj2} a obtenu une loi de Bohr-Sommerfeld  
complexe qui montre que les valeurs propres sont situ\'es sur un
r\'eseau d\'eform\'e. D'autre part, avec M.~Hager \cite{HaSj}
nous avons montr\'e dans toute dimension que l'asymptotique de 
Weyl est valable avec probabilit\'e proche de 1 
pour des petites perturbations al\'eatoires de l'op\'erateur.
Dans les deux cas, le valeurs propres sont distribu\'ees (\`a des
petites corrections pr\`es) selon des densit\'es lisses, et ici nous
montrons que ces densit\'es sont en g\'en\'eral diff\'erentes. 
\end{abstract}


\section{Introduction}
\label{al}
\setcounter{equation}{0}

In the classical paper by J.J.~Duistermaat and L.~H\"ormander
\cite{DuHo}, one very interesting application is about
(pseudo)differential operators with principal symbol $p$ such that the
Poisson bracket $\{ p,\overline{p}\}$ vanishes on the zero set
$p^{-1}(0)$ and the differentials of the real and imaginary part of $p$
are independent there, so that the zero set is a codimension 2
sub-manifold of the cotangent space. The authors gave interesting
existence results under non-compactness assumptions on the
bicharacteristic foliation. In my thesis under the direction of
L.~H\"ormander my task was to study the case when $\{
p,\overline{p}\}\ne 0$ on the zero-set of the symbol and in a
subsequent paper with Duistermaat \cite{DuSj} we 
introduced an studied certain microlocal projections onto
the kernel and the co-kernel of the operator. The full history of this
subject can be traced back to the famous counterexample of
Hans Lewy to local solvability and subsequent works by H\"ormander and
others and there is also quite a rich recent history.  

There has been a renewed interest in non-self-adjoint
operators and the related notion of pseudospectrum, promoted by
L.N.~Trefethen, E.B.~Davies, M.~Zworski and others. Again the
Poisson bracket $i^{-1}\{ p,\overline{p}\}$ plays an important role as a
source of pseudospectral behaviour, including spectral instability. (We
observe here that the above Poisson bracket is equal to the principal
symbol of the commutator of the corresponding (pseudo)differential
operator and its adjoint.) We refer to the surveys \cite{Sj1, Sj2} 
where further references can be found.

Possibly, as a reaction to these developments, the author participated 
in two projects: 
\begin{itemize}
\item With A.~Melin \cite{MeSj2} we discovered for a a fairly wide and
  stable class of non-self-adjoint semi-classical pseudodifferential
  operators with analytic symbols that the individual eigenvalues in
  certain regions can be determined by a Bohr-Sommerfeld quantization
  rule defined in terms of certain complex Lagrangian tori (close to the real
  domain). The underlying idea is here to change the Hilbert space
  norm by means of exponential weights in such a way that the operator
  becomes (more) normal.
\item M.~Hager \cite{Ha} considered certain non-self-adjoint
  $h$-pseudodifferential operators in dimension 1 with small
  multiplicative random perturbations and showed that with probability
  tending to 1 when $h\to 0$, the eigenvalues distribute according to
  the classical Weyl law, well-known in the context of self-adjoint
  operators since almost a century. The same type of result was
  subsequently obtained in any dimension by Hager and the author
  \cite{HaSj} for a certain class of non-multiplicative random perturbations
  and recently also for multiplicative random perturbations in any
  dimension by the author \cite{Sj3}.
\end{itemize}
In the present note we shall compare the resulting distributions of
eigenvalues, in the case when the unperturbed operator satisfies the
assumptions of \cite{MeSj2}, and we shall see that they are in general
different. This means that the random perturbations will change
radically the asymptotic distribution of eigenvalues. 
The intuitive explanation of this
phenomenon is that the result of \cite{MeSj2} depends on the geometry
in the complex domain, while the random perturbation destroys
analyticity and hence the eigenvalue distribution should be given in
terms of the real phase space, where the Weyl law is the natural
candidate.

We next describe the main result of \cite{MeSj2}.
Let $p(x, \xi )$ be bounded and holomorphic in a tubular 
neighborhood of ${\bf R}^4$ in
${\bf C}^4 = {\bf C}^2_x\times{\bf C}^2_\xi $.
(The assumptions near $\infty $ can be varied in many ways and we can
let $p$ belong to some more general symbol space as long as we have
the appropriate form of ellipticity near infinity, cf \no{int.2}
below.) 
Assume that
\ekv{int.1}{{\bf R}^4 \cap p^{-1}(0)\ne \emptyset \mbox{ is
    connected,}}
for simplicity. Also assume that
\ekv{int.2}{\mbox{on } {\bf R}^4\mbox{ we have }|p(x, \xi )| \ge
1/C,\mbox{ for }|(x, \xi )| \ge  C,}
for some $C > 0$,
\ekv{int.3}{d\Re p(x, \xi ), d \Im p(x, \xi )
\mbox{ are linearly independent for all } (x, \xi ) \in p^{-1}(0) \cap
{\bf R}^4.}
It follows that $p^{-1}(0)\cap {\bf R}^4$ is a compact (2-dimensional) 
surface. 

\par Also assume that 
\ekv{int.4}{|\{ \Re p, \Im p\} |\mbox{ is
sufficiently small on } p^{-1}(0) \cap {\bf R}^4.}  
By ``sufficiently
small'', we mean that $| \{ \Re p,\Im p\} | <\delta $ for some $\delta
>0$ that will depend on all constants (implicit or explicit) that are
required to express the other conditions above uniformly.

In \cite{MeSj2} we showed that $p^{-1}(z)\cap {\bf R}^4$ is a real
torus for $z\in \mathrm{neigh\,}(0,{\bf C})$ (ie some neighborhood of $0$ in
${\bf C}$) and that there exists a smooth
$2$-dimensional torus $\Gamma (z)\subset p^{-1}(z)\cap{\bf C}^4$,
close to $p^{-1}(z)\cap{\bf R}^4$ such that ${{\sigma }_\vert}_{\Gamma
(z)}=0$ and $I_j(z)\in {\bf R}$, $j=1,2$. Here $I_j(z):=\int_{\gamma
_j(z)} \xi \cdot dx$ are the actions along two fundamental cycles
$\gamma _1(z), \gamma _2(z)\subset \Gamma (z)$ and $\sigma =\sum_{1}^2
d\xi _j\wedge dx_j$ is the complex symplectic $(2,0)$ form. Moreover, 
$\Gamma (z)$, $I_j(z)$ depend smoothly on $z\in \mathrm{neigh\,}(0)$.

\par The main result of \cite{MeSj2}, valid under slightly more
general assumptions than the ones above, is then
\begin{theo}\label{int1}
Under the above assumptions, there exist a neighborhood $V$ of $0\in
{\bf C}$, $\theta_0 \in ( \frac{1}{2}
{\bf Z})^2$, $\theta _j\in C^\infty (V;{\bf R}^2)$ and $\theta (z; h) \sim
\theta_ 0 + \theta_ 1(z)h + \theta_ 2(z)h^2 + ..$ in $C^\infty (
V;{\bf R}^2)$, such that for $z\in V$
and for $h > 0$ sufficiently small, 
$z$ is an eigenvalue of $P=p^w(x,hD_x)$ iff
$$\frac{(I_1(z),I_2(z))}{2\pi h}
= k-\theta (z; h),\mbox{ for some }k \in {\bf Z}^2.\quad (BS)$$
Here $p^w(x,hD)$ denotes the Weyl quantization of the symbol $p(x,h\xi )$.
\end{theo}

\par Let us also assume that 
\ekv{int.4.5}{
\mbox{the map }z\mapsto
I(z):=(I_1(z), I_2(z))\mbox{ is a diffeomorphism from }V\mbox{ to }
I(V).} 
This assumption is satisfied if we strengthen (\ref{int.4}) by
assuming that $|\{\Re p,\Im p\}|$ is sufficiently small on $p^{-1}(z)$
 for all $z\in
\mathrm{neigh\,}(0,{\bf C})$ and choose $V$ small enough. The
eigenvalues near $0$ will then form a distorted lattice and we
introduce the leading spectral density function $0<\omega (z)\in 
C^\infty (V)$ by 
\ekv{int.5}
{
dI_1(z)\wedge dI_2(z)=\pm \omega (z) d\Re z\wedge d\Im z,
} 
where the sign is chosen so that $\omega $ becomes positive. Then from
Theorem \ref{int1} it follows that for every $W\Subset V$ with smooth
boundary, the number of eigenvalues in $W$ satisfies
\ekv{int.6}
{
N(W;h)=\frac{1}{(2\pi h)^2}(\int_W \omega (z)L(dz)+o(1)),\ h\to 0.
}
Here $L(dz)=d\Re zd\Im z$ denotes the Lebesgue measure.

Now we turn to the results in \cite{Ha, HaSj, Sj3}. Again the
unperturbed operator is of the form $P=p^w(x,hD_x)$ where the
complex-valued smooth symbol should belong to a suitable symbol class
and satisfy an ellipticity condition at infinity which guarantees that
the spectrum of $P$ in a given open set $\Omega \Subset {\bf C}$ is
discrete. The perturbed operator is of the form $P_\delta =P+\delta
Q_\omega $, where the parameter $\delta $ is small, say bounded from
above by some positive power of $h$ and from below by $e^{-h^{-\alpha
  }}$ for some suitable value $\alpha \in ]0,1]$. Under some additional
assumptions on the type of random perturbation and about
non-constancy of the symbol $p$, it is showed in the cited works that
with a probability that tends to 1 when $h\to 0$, the number of
eigenvalues of $P_\delta $ in $W\Subset \Omega $ obeys 

\ekv{int.6.2}
{
N_\delta (W;h)=\frac{1}{(2\pi h)^n}(\mathrm{vol}\,(p^{-1}(W))+o(1))
}
uniformly for $W$ in a class of subsets of $\Omega $ 
with uniformly smooth boundary. (In the case of multiplicative
perturbations, an additional symmetry assumption on the symbol is
imposed which cannot be completely eliminated.)

\par Notice that this result can be formulated as in (\ref{int.6})
with the density $\omega $ replaced by the Weyl density $w(z)L(dz)$,
defined to be the direct image of the symplectic volume element under 
the map $p$, so that 
\ekv{int.7}
{
\int f(z)w(z)L(dz)=\iint f(p(x,\xi ))dxd\xi ,\ f\in C_0^\infty (V).
}

\par 
In the two-dimensional case there are situations (for instance in the
case of the symbol $p(x,\xi )=\frac{1}{2}((x_1^2+\xi _1^2)+i(x_2^2+\xi
_2^2))-\mathrm{const.}$ and small perturbations of that symbol) where
Theorem \ref{int1} applies to $P $ and the results of
\cite{HaSj, Sj} apply to small random perturbations, and it is then of
interest to compare the spectral densities.
We shall see that $\omega (z)=w(z)$ in the integrable case, when
$\{ \Re p,\Im p\}\equiv 0$ but that these quantities are 
different  in general.
\begin{theo}\label{int2}
Under the assumptions (\ref{int.1})--(\ref{int.4.5}) 
we have generically that 
$w\not\equiv \omega $.
\end{theo}
In other words, if $w\equiv \omega $, 
then there are arbitrarily small perturbations of $P$
within the class of operators as in the theorem, 
for which $w\not\equiv \omega $.
\section{The integrable case}\label{ic}
\setcounter{equation}{0}
In this section, we strengthen the assumption \no{int.4} to 
\ekv{ic.1}
{
\{\Re p,\Im p\} \equiv 0.
}
It is then well-known by the Liouville-Mineur-Arnold theorem (see \cite{Vu}) that
there exists a real symplectic diffeomorphism 
$\kappa :\mathrm{neigh\,}(\eta =0,T^*{\bf T}^2)\to \mathrm{neigh
\,}(p^{-1}(0)\cap {\bf R}^4,{\bf R}^4)$, (i.e. from a neighborhood of $\{ \eta
=0\}$ in $T^*{\bf T}^2$ to a neighborhood of $p^{-1}(0)\cap {\bf R}^4$
in ${\bf R}^4$) such that 
\ekv{ic.2}
{
p\circ \kappa =\widetilde{p}(\eta )
}
is independent of $y$, where ${\bf T}^2=({\bf R}/2\pi {\bf Z})^2$ and 
$T^*{\bf T}^2\simeq {\bf T}_y^2\times {\bf R}_\eta ^2$.

\par In this case $\Gamma (z)$ is simply the real Lagrangian torus
$p^{-1}(z)\cap {\bf R}^4$ and 
\ekv{ic.3}
{
I_j(z)=2\pi \eta _j+I_j(0),\quad \widetilde{p}(\eta )=z.
}
It follows that up to composition with
$\kappa $, we get the same quantities $\omega (z)$, $w(z)$, if we do
the computation directly on $T^*{\bf T}^2$ for $\widetilde{p}(\eta )$,
and we restrict the attention to that case
and drop the tilde.

\par From \no{ic.3} and \no{int.5} we get 
\ekv{ic.4}
{
\frac{\omega (z)}{(2\pi )^2}=\left| \det \frac{\partial (\eta _1,\eta
    _2)}{\partial (\Re p,\Im p)}\right| ,\quad p(\eta )=z.
}

\par From \no{int.7}, we get for $f\in C_0^\infty (V)$:
\begin{eqnarray*}
&\displaystyle\int f(z) w(z)L(dz)=\iint f(p(y,\eta ))dyd\eta& \\
&\displaystyle=(2\pi )^2\int f(p(\eta ))d\eta =(2\pi )^2\int
f(z)\left| \det \frac{\partial (\eta _1,\eta
    _2)}{\partial (\Re p,\Im p)}\right|L(dz),&
\end{eqnarray*}
which shows that $w(z)$ also satisfies \no{ic.4}, so
\ekv{ic.5}
{
w(z)=\omega (z),\ z\in V,
} 
in the completely integrable case \no{ic.1}.  

\section{The general case}\label{gc}
\setcounter{equation}{0}

\par In this section we shall prove Theorem \ref{int2} by means of
calculations similar to the ones in \cite{MeSj1}. Let $p$ satisfy the 
assumptions (\ref{int.1})--(\ref{int.4.5}). Let $G_t(x,\xi )$ for $t\in
\mathrm{neigh\,}(0,{\bf R})$ be a smooth family of functions that are
holomorphic and uniformly bounded in a fixed tubular neighborhood of
${\bf R}^4$. Possibly after decreasing the neighborhood of $t=0$ we
get a smooth family of canonical transformations $\kappa _t$
from a fixed tubular neighborhood of ${\bf R}^4$ onto a neighborhood
of ${\bf R}^4$, by solving the Cauchy problem 

\ekv{gc.3}{\frac{d}{dt}\kappa _t(\rho )=(\kappa _t)_*(\widehat{iH_{G_t}})(\rho ),\quad \kappa
  _0(\rho )=\rho ,}
where $H_{G_t}=\frac{\partial G_t}{\partial \xi }\frac{\partial
}{\partial x}-\frac{\partial G_t}{\partial x }\frac{\partial
}{\partial \xi }$ is the holomorphic Hamilton field (of type (1,0))
and we identify $iH_{G_t}$ with the corresponding real vector field 
$\widehat{iH_{G_t}}:=iH_{G_t}+\overline{iH_{G_t}}$. 

\par Put $p_t=p\circ \kappa _t$. Then (possibly after further shrinking
the neighborhood of $t=0$) $p_t$ will satisfy the assumptions
(\ref{int.1})--(\ref{int.4.5}) and since $\kappa _t$ are complex canonical
transformations, we also know that 
\ekv{gc.4}
{
\omega _t = \omega \hbox{ is independent of }t.
}
In order to prove Theorem \ref{int2}, it suffices to show
\begin{theo}\label{gc1}
For every neighborhood $V$ of $0\in {\bf C}$, we can find a family
$G_t$ as above, such that every neighborhood of $t=0$ will contain a
$t$ for which $w_t\not\equiv w$ in $V$. Here $w_t$ denotes the Weyl
density of $p_t$, defined as in (\ref{int.7}).
\end{theo}

\begin{remark}\label{gc2}
{\rm Actually, we shall prove the theorem in all dimensions (replacing 2 by
any $0<n\in {\bf N}$) for any $p\not\equiv 0$ that is bounded and
holomorphic in a tubular neighborhood of ${\bf R}^{2n}$ in ${\bf
  C}^{2n}$ that satisfies (\ref{int.2}) and for which
$p^{-1}(0)\cap {\bf R}^{2n}\ne \emptyset $. 
Let $w_t(z)L(dz)$ be the measure
defined as in (\ref{int.7}) with $p=p_0$ replaced by $p_t=p\circ
\kappa _t$. }
\end{remark}
\begin{proof}
 For $f\in C_0^\infty (V;{\bf R})$ we get,
\begin{eqnarray*}
&&\int f(z)\frac{\partial w_t(z)}{\partial t}L(dz)=\frac{d}{dt} \int
f(z)w_t(z)L(dz)=\\
&&\frac{d}{dt}\iint f(p_t(x,\xi ))dxd\xi =\iint (\frac{\partial
  f}{\partial z}(p_t)\frac{\partial p_t}{\partial t}+\frac{\partial
  f}{\partial \overline{z}}(p_t)\frac{\partial \overline{p}_t}{\partial t})dxd\xi .
\end{eqnarray*}
Here, we have 
$$
\frac{\partial p_t}{\partial t}=iH_{G_t}p_t,
$$
and using that $f$ is real,
\begin{eqnarray}\label{gc.5}
\int f(z) \frac{\partial w_t(z)}{\partial t}L(dz)&=&2\Re (i\iint
\frac{\partial f}{\partial z}(p_t)H_{G_t}p_t dxd\xi)\\ 
&=&-2\Re (i \iint \frac{\partial f}{\partial z}(p_t)H_{p_t}(G_t)dxd\xi
)
\nonumber \\
&=& 2\Re (i \iint H_{p_t}\left(\frac{\partial f}{\partial z}(p_t) 
\right) G_t
dxd\xi )\nonumber\\
&=& 2\Re (i\iint (\frac{\partial ^2f}{\partial z^2}H_{p_t}(p_t)+ 
\frac{\partial ^2f}{\partial \overline{z}\partial z}H_{p_t}(\overline{p}_t)G_t
dxd\xi )\nonumber\\
&=& 2\Re (\iint \frac{\partial ^2f}{\partial \overline{z}\partial
  z}(p_t) i\{ p_t,\overline{p}_t\} G_tdxd\xi )\nonumber\\
&=& \iint (\Delta f)(p_t)\{ \Re p_t,\Im p_t\} \Re G_t dxd\xi .\nonumber
\end{eqnarray}
If $\{\Re p,\Im p\} =\frac{i}{2}\{ p,\overline{p}\}$ does not vanish
identically, there are points arbitrarily close to $p^{-1}(0)$ where it
does not vanish and we can choose $f\in C_0^\infty (V;{\bf R})$ (where
$V$ is any fixed neighborhood of $0\in {\bf C}$) such that $(\Delta
f)(p)\{\Re p,\Im p\}$ does not vanish identically. We can then choose
$G=G_0$ independent of $t$ with the properties above so that 
$$\int f(z)\left( \frac{\partial }{\partial t}\right)_{t=0}w_t(z)L(dz)
=\iint (\Delta f)(p)\{ \Re p,\Im p\} \Re G dxd\xi \ne 0.
$$
We get the conclusion of Theorem \ref{gc1} in this case.

\par If $\{ \Re p,\Im p\} \equiv 0$, we choose $G$ real and
independent of $t$ in (\ref{gc.5}) and differentiate that identity
once with respect to $t$ at $t=0$ to get:
\begin{eqnarray*}
\int f(z)\left(\frac{\partial ^2w_t}{\partial
    t^2}\right)_{t=0}L(dz)&=&
\iint (\Delta f)(p)\left(\frac{\partial }{\partial
    t}\right)_{t=0}(\frac{i}{2}
\{ p_t,\overline{p_t}\} )G dxd\xi \\
&=& \iint (\Delta f)(p)\frac{i}{2}(\{ iH_Gp,\overline{p}\} +\{
p,\overline{iH_Gp}\})Gdxd\xi \\
&=& -\frac{1}{2}\iint (\Delta
f)(p)(H_{\overline{p}}H_pG+H_pH_{\overline{p}}G)Gdxd\xi .
\end{eqnarray*}
Here we integrate by parts and use that $H_pp=0$,
$H_{\overline{p}}p=0$, to get
$$
\int f(z)\left(\frac{\partial ^2w_t}{\partial t^2}\right)_{t=0}L(dz)
=
\iint (\Delta f)(p)|H_pG|^2 dxd\xi .
$$ 
Again we see that we can find $f\in C_0^\infty (V;{\bf R})$ and $G=G_0$
as above, so that the last integral is $\ne 0$. The conclusion in the
theorem follows in this case also.
\end{proof}

\end{document}